\documentclass[11pt]{article}
\usepackage{latexsym}
	\usepackage{amscd}
	\usepackage{amssymb}
	\usepackage{amsmath}
	\usepackage{amsthm}
	\usepackage{enumerate}
	\usepackage{verbatim}
\def\co{\colon\thinspace}

\newtheorem{thm}{Theorem}[section]

\newtheorem{cor}[thm]{Corollary}
\newtheorem{quest}[thm]{Question}
\newtheorem{lem}[thm]{Lemma}

\newtheorem{Example}[thm]{Example}
\newenvironment{ex}{\begin{Example}\rm}{\end{Example}}
\newtheorem{Counterexample}[thm]{Counterexample}
\newenvironment{coex}{\begin{Counterexample}\rm}{\end{Counterexample}}
\newtheorem{remark}[thm]{Remark}
\newenvironment{rmk}{\begin{remark}\rm}{\end{remark}}
\newtheorem{Fact}[thm]{Fact}

\newtheorem{Nothing}[thm]{$\!\!\!$}

\newcommand{\be}{\begin{equation}}
\newcommand{\ee}{ \end{equation}} 
\newcommand{\ba}{\begin{eqnarray}}
\newcommand{\ea}{\end{eqnarray}}
\newcommand{\ban}{\begin{eqnarray*}}
\newcommand{\ean}{\end{eqnarray*}}

\newcommand{\lp}{\langle}
\newcommand{\rp}{\rangle}

\newcommand{\Ric}{\mbox{Ric}}

\begin{document}
\abovedisplayskip=6pt plus3pt minus3pt
\belowdisplayskip=6pt plus3pt minus3pt
\title{\bf Metrics of positive Ricci curvature on vector bundles 
over nilmanifolds\rm}
\author{Igor Belegradek \and 
Guofang Wei\thanks {Partially supported by NSF Grant \# DMS-9971833.}}
\date{}
\maketitle
\begin{abstract}
We construct metrics of positive Ricci curvature on some vector
bundles over tori (or more generally, over nilmanifolds). This
gives rise to the first examples of manifolds with positive Ricci curvature 
which are homotopy equivalent but not homeomorphic to manifolds
of nonnegative sectional curvature.
\end{abstract}

\section{Introduction}
According to the soul
theorem of J.~Cheeger and D.~Gromoll, a complete open manifold of
nonnegative sectional curvature is diffeomorphic to the total space of the
normal bundle of a compact totally geodesic submanifold.
Furthermore, a finite cover of any compact manifold with $sec\ge 0$
is diffeomorphic to the product $C\times T$ where $T$ is a torus,
and $C$ is simply-connected with $sec\ge 0$. (The same is actually
true when $sec$ is replaced by $Ric$.)

It was shown in~\cite{OW, BK2, BK3} that a majority of vector bundles over
$C\times T$ admit no metric with $sec\ge 0$ provided $\dim(T)$ is
sufficiently large.
For example, any vector bundle over a torus whose total space admits a 
metric with $sec\ge 0$ becomes trivial in a finite cover, and
there are only finitely many such bundles in each rank.

It is then a natural question which of the bundles can carry metrics with
$Ric\ge 0$ or $Ric>0$. J.~Nash and L.~Berard-Bergery~\cite{Nash, BerB}
constructed metrics with $Ric>0$ on every vector bundle of rank $\ge 2$
over a compact manifold $C$ with $Ric(C)>0$. 
Their method generally fails if $C$ only
has $Ric(C)\ge 0$. In fact, the metrics produced in~\cite{Nash, BerB}
have the property that the bundle projection is a Riemannian submersion, 
while according to~\cite[2.2]{BK2}, any vector bundle over a torus
with $\Ric\ge 0$ on the total space becomes trivial in a finite cover
provided the bundle projection is a Riemannian submersion.

Recall that a {\it nilmanifold}
is the quotient of a simply-connected nilpotent group by a
discrete subgroup. Any nilmanifold is diffeomorphic
to the product of a compact nilmanifold and a Euclidean space~\cite{Mal}.
Compact nilmanifolds can be described as iterated principal 
circle bundles~\cite{OT}. 
Here is our main result 
(see also~\ref{main thm tech} for
a more technical version).

\begin{thm}\label{main thm}
Let $L$ be a complex line bundle over a  nilmanifold.
Then for all sufficiently large $k$, the total space of the Whitney sum
of $k$ copies of $L$ admits a complete Riemannian metric of positive
Ricci curvature.
\end{thm}

Note that if the total spaces $E(\xi)$, $E(\eta)$ of vector bundles
of $\xi$, $\eta$ over a smooth aspherical manifold $M$ admit complete
metrics with $Ric\ge 0$, then so does the total space of 
$\xi\oplus \eta\oplus TM$ because $E(\xi\oplus \eta\oplus TM)$
is the covering of $E(\xi)\times E(\eta)$ corresponding to
the subgroup $d_*\pi_1(M)$ where $d\co M\to M\times M$
is the diagonal embedding. Since complete nilmanifolds are aspherical
and parallelizable, Theorem~\ref{main thm} generalizes,
up to taking sum with some trivial bundle, to vector bundles 
with tori as structure groups (the structure group of a vector bundle
is a $k$-torus iff the bundle is a Whitney sum of $k$ complex line bundles
and a trivial bundle). 
Combining~\ref{main thm tech} with~\cite[1.2]{BK2}, 
we deduce the following.

\begin{cor}\label{prod with comp}
Let $C$ be a point or a compact manifold with $Ric\ge 0$, and
let $T$ be a torus of dimension $\ge 4$. Then for all sufficiently large
$k$, there is an infinite sequence of rank $k$ vector bundles 
over $C\times T$ whose total spaces are pairwise nonhomeomorphic,
admit complete Riemannian metrics with $Ric>0$,
and are not homeomorphic to complete Riemannian  manifolds with $sec\ge 0$.
\end{cor}

Corollary~\ref{prod with comp} gives a new type of examples
of manifolds with $Ric>0$. These manifolds are topologically
similar but not homeomorphic to manifolds with $sec\ge 0$.
We refer to~\cite{Grl} for the survey of known examples.
Now the following question is natural.

\begin{quest}\label{quest}
Do most vector bundles over 
tori (or more generally compact manifolds with
$Ric\ge 0$ or $sec\ge 0$) admit metrics with $Ric\ge 0$, or
$Ric>0$?
\end{quest}

Note that by M.~Anderson's first Betti number estimate~\cite{And}, no
vector bundle of rank $\le 2$ over a torus admits a metric with
$Ric>0$. This however seems to be the only known obstruction.
Theorem~\ref{main thm} also gives metrics of $Ric>0$ on
some nontrivial vector bundles over many nilmanifolds as follows. 

\begin{cor}\label{cor nil}
Let $N$ be a nilmanifold such that there exists a class 
in $H^2(N,\mathbb Z)$ with nontorsion cup-square. 
Then for all sufficiently large
$k$, there is an infinite sequence of rank $k$ vector bundles 
over $N$ whose total spaces are pairwise nonhomeomorphic and
admit complete Riemannian metrics with $Ric>0$.
\end{cor}

The existence of a $2$-class with a nontorsion cup-square
is equivalent, by~\ref{lem ina app}, 
to the existence of a complex line bundle
$L$ over $N$ such that for any $k$, the Whitney sum of
$k$ copies of $L$ is nontrivial as a real vector bundle.
In~\ref{five-nil} we give an example of a compact nilmanifold
of dimension $5$ for which such an $L$ does not exist.

Metrics of positive Ricci curvature on vector bundles over nilmanifolds
are interesting in their own right. Apparently, besides
Anderson's growth estimate~\cite{And}, no obstructions
are known to the existence of such metrics.

Here is the idea of the proof of Theorem~\ref{main thm}. The second author
showed in~\cite{W} that if $G$ is a simply-connected nilpotent Lie 
group, then $G\times\mathbb R^k$ has a warped product
metric with $Ric>0$ such that $G$ lies in the isometry group.
In this paper we start with a similar metric on $G\times\mathbb C^k$
where $G$ has a central one-parametric subgroup $Q$. 
Let $Q$ act on $G$ by left multiplication,
fix an epimorphism $Q\to S^1$, and a 
diagonal isometric action of $S^1$ on $\mathbb C^k$.
We then show that the Riemannian submersion metric on
$(G\times\mathbb C^k)/Q$ has $Ric>0$ for all large enough $k$.
(Note that, generally, Ricci curvature may decrease under Riemannian
submersions, namely, there is an error term in the O'Neill submersion 
formula for Ricci curvature. The error term always vanishes 
for submersions with minimal fibers, however, in our case 
the fibers clearly have different volume on the 
$\mathbb C^k$-factor so they are not minimal.  
We prove that $(G\times\mathbb C^k)/Q$ has $Ric>0$
by being able to control the error term when $k$ is sufficiently large.)
The obvious $G$-action on  $(G\times\mathbb C^k)/Q$ is isometric so
the quotients of  $(G\times\mathbb C^k)/Q$ by discrete subgroups
of $G$ also carry metrics with $\Ric>0$.
The bundles described in~\ref{main thm} can be obtained as  
$(P\times\mathbb C^k)/Q$ where $P$ is a principal circle bundle
associated with $L$.
Finally,~\ref{main thm} follows from 
the fact that a principal circle
bundle over a nilmanifold is again a nilmanifold.

The structure of the paper is as follows.
In sections~\ref{sec: def metric} and~\ref{sec: set up}
we define the metric, and choose a coordinate system
to compute the Ricci curvature. This computation is done
in sections~\ref{sec: error} and~\ref{sec: mixed terms}.
Applications are discussed in 
section~\ref{sec: metrics on bundle}.
In Appendix we prove some facts 
on the existence of line bundles over nilmanifolds.

The first author is grateful to Vitali Kapovitch
for motivating discussions.

\section{Defining the metric}
\label{sec: def metric}
Let $G$ be a simply-connected nilpotent Lie group with the Lie
algebra $g$. 
Fix an arbitrary central element $X_1\in g$ and
use the ascending central series for $g$ to
choose a basis $X_1,\cdots, X_n$ for $g$ such that
$[X, X_i]\in\mbox{span}\{X_1, \cdots, X_{i-1}\}$ for all $X\in g$. 
In particular,
\begin{equation}\label{xbracket}
[X_1, X_i] = 0 \ \mbox{for all} \ i = 1, \cdots, n.
\end{equation}

Equip $G\times {\mathbb C}^k$ with the warped product metric 
constructed in~\cite{W}. Namely, let
\[ g = g_r + dr^2 + f^2 (r) ds^2_{2k-1}  \]
be the metric on $G\times {\mathbb C}^k$, where 
$f(r) = r(1+r^2)^{-1/4}$, $ds^2_{2k-1}$ is the canonical metric 
on the sphere $S^{2k-1}$ 
(the one induced by the inclusion of the unit sphere into $\mathbb  C^k$), 
and $g_r$ is a special family of almost 
flat left invariant metrics on $G$ defined as follows. 

Let $Y_i = X_i/h_i (r)$, where 
\begin{equation} \label{const-alphai}
h_i (r) = (1+r^2)^{-\alpha_i},
\end{equation}
and $\alpha_i = 2^{n-i+1} + 2^{-1-i} -\frac 12$. 
(This value of $\alpha_i$ is obtained by
choosing $\alpha_n = 2^{-1-n} + \frac 32$ in~\cite{W}.) 
Then let $g_r$ be the left-invariant metric on $G$ such that 
$\{ Y_1, \cdots, Y_n \}$ form an orthonormal basis. 
As we explain in~\ref{ric-yi-yj}, the Ricci tensor of 
the metric $g_r$ satisfies
\begin{equation} \label{const-c} 
|\Ric_G (Y_i, Y_j) | \leq c (1+r^2)^{-1},
\end{equation}
where $c$ is a constant depending on $n$ and the structure constants.
Clearly, $g$ is a complete metric on $G\times {\mathbb  C}^k$ whose 
isometry group contains $G$. 
By~\cite{W}, $\Ric (g)>0$ for all sufficiently large $k$.

For a nonzero integer $m$, let 
$\varphi^m\co\mathbb R\to U(1)$ be the homomorphism
given by $\varphi^m(q)=e^{imq}$. 
Define a sequence of isometric 
$\mathbb R$-actions $\Phi_m$ on $G\times\mathbb C^k$ by
\begin{equation}    
\label{action}
\Phi_m(q)(x,u)=
(\mathrm{exp}(qX_1)x,\varphi^m(q)u) 
\ \mbox{for any} \ x\in G, u\in\mathbb C^k.
\end{equation}
Equip $(G\times{\mathbb C}^k)/\Phi_m(\mathbb R)$ 
with the Riemannian submersion metric.
Note that the isometric 
$G$-action on $G\times\mathbb C^k$
given by $g(h,z)=(gh,z)$ induces an isometric
$G$-action on the quotient $(G\times{\mathbb C}^k)/\Phi_m(\mathbb R)$;
the kernel of the action is exactly
$K_m=\mathrm{exp}(\frac{2\pi}{m}\mathbb Z X_1)$ so that
$G/K_m$ acts effectively.

\section{Setting up the computation}
\label{sec: set up}
Let $\Ric$ denote the Ricci curvature on the total space 
$G\times {\mathbb C}^k$ of the submersion
and $\check{\Ric}$ the Ricci curvature on the base 
$(G \times {\mathbb  C}^k)/\Phi_m(\mathbb R)$. 
Then for any vector fields $X,Y$ on the base, the Ricci curvatures 
are related  by the following (see~\cite[p.244]{Bes})
\begin{equation}\label{ric}
\check{\Ric} (X,Y) = 
\Ric (X,Y) + 2 \langle A_X, A_Y \rp + \langle TX, TY \rp -\frac 12 
(\langle \nabla_X S, Y \rp + \langle \nabla_Y S, X \rp ),
\end{equation}
where $A,T$ are the O'Neill tensors, and $S$ is the mean curvature 
of the fiber (for precise definition, see, 
e.g.~\cite[p.239]{Bes}). In our case 
$\langle A_X, A_Y \rp = \langle A_X W, A_Y W \rp, \ 
\langle TX, TY \rp = \langle T_W X, T_W Y \rp$
where $W$ is defined in (\ref{w}). 

We want to show that for all sufficiently large $k$, 
the Riemannian submersion metric on 
$(G\times{\mathbb C}^k)/\Phi_m(\mathbb R)$ has positive 
Ricci curvature. 
To prove this, it suffices to control the error term 
$\check{\Ric}(X,Y)-{\Ric}(X,Y)$ in (\ref{ric})
coming from the 
mean curvature tensor of the fiber of the $\mathbb R$-action.

For the rest of the computation, fix a point $(g_*, r_*, s_*)$
of $G\times\mathbb C^k$
and show that the horizontal space tangent to the point
has all Ricci curvatures positive.

At the point $s_*$ of the sphere $(S^{2k-1}, ds^2_{2k-1})$ 
choose an orthonormal frame
$V_1, \cdots, V_{2k-1}$ with $V_1$ tangent to the orbit of the 
action $\Phi_m$ restricted to $S^{2k-1}$ such that 
\begin{equation}\label{vbracket}
{\nabla_{V_i}V_j}\vert_{s_*} = 0 \ \mbox{for all} \ i,j = 1, \cdots, 2k-1.
\end{equation}
%

Let $U_j = V_j/f(r)$. Then $\partial_r = \frac{\partial}{\partial r}, U_1, 
\cdots U_{2k-1}, Y_1, \cdots, Y_n$ form an orthonormal basis of $g$ on 
$G\times {\mathbb  C}^k$ in a neighborhood of $(g_*, r_*, s_*)$.
The vector field $X_1 + m V_1$ is tangent to the orbits of
the $\mathbb R$-action $\Phi_m$ on $G\times {\mathbb  C}^k$.
Let 
\begin{equation}\label{w}
W = (X_1 + m V_1)/({h_1}^2 + m^2 f^2)^{1/2} = 
(h_1 Y_1 + m f U_1)/({h_1}^2 + m^2f^2)^{1/2}
\end{equation}
be its unit vector. 
Then $\partial_r, U_2, \cdots U_{2k-1}, Y_2, \cdots, Y_n, W'$ 
form an orthonormal basis 
of $(N \times {\mathbb  C}^k)/\Phi_m(\mathbb R)$, where 
\begin{equation}\label{w'}
W' = (m f Y_1-h_1 U_1)/({h_1}^2 + m^2 f^2)^{1/2}.
\end{equation}

The mean curvature vector $S$ of the fiber is $S = (\nabla_W W)^H$, the 
horizontal projection of $\nabla_W W$.
Using the Koszul's formula for the Riemannian connection 
(the one used to prove its existence) and (\ref{xbracket}), 
(\ref{vbracket}), we have 
\[  \nabla_{Y_1} Y_1 = -\frac{h_1'}{h_1} \partial_r,
\ \ \  \nabla_{U_1} U_1 = 
-\frac{f'}{f} \partial_r, \ \ \ 
\nabla_{Y_1} U_1 = \nabla_{U_1} Y_1 = 0.\]
Since any function of $r$ has zero derivatives in 
the directions $Y_1$, $U_1$, the Koszul's formula implies
\begin{equation} \label{s} 
S= - \frac{h_1 h_1' + m^2 f f'}{{h_1}^2 + m^2 f^2} \partial_r.
\end{equation}

The Ricci curvature on the total space $G\times\mathbb C^k$
are given by~\cite{W, W2} and~\ref{ric-yi-yj}.
\begin{eqnarray*}
\label{mixed-total}
 & \Ric(\partial_r, U_{j}) = 0  & (1 \leq j \leq 2k-1), \\
 & \Ric(Y_{i}, \partial_r) = \Ric(Y_{i}, U_{j}) = 0  
 & (1 \leq i \leq n,\ 1 \leq 2k-1), \\
 & |\Ric(Y_{i}, Y_{j})| \leq   c (1+r^2)^{-1},  & 
(i \not= j, 1 \leq i, j \leq n),
\end{eqnarray*}
where $c$ is a constant depending only on $n$ and the structure constants.
\begin{eqnarray}
\Ric(\partial_r, \partial_r) & = 
& -\sum_{i=1}^{n}\frac{h''_{i}(r)}{h_{i}(r)} 
-(2k-1)\frac{f''(r)}{f(r)}, \nonumber \\
& = & \left\{- \sum_{i=1}^{n} 2\alpha_{i}[(2\alpha_{i} + 1)r^{2} - 1] 
  + (2k-1)\frac{r^{2} + 6}{4} \right\}/(1 + r^{2})^{2}.  \label{p-r}\\
\vspace*{.1in}
\Ric(Y_{i}, Y_{i}) & = & -\frac{h''_{i}(r)}{h_{i}(r)} - (2k-1)\frac{h'_{i}(r)
f'(r)}{h_{i}(r)f(r)} + Ric_{G} (Y_{i}) -\sum_{i\not= j}\frac{h'_{i}(r)h'_{j
}(r)}{h_{i}(r)h_{j}(r)} \nonumber \\
& \geq & \{ -2\alpha_{i}[(2\alpha_{i} + 1
)r^{2} - 1] + (2k-1)\alpha_{i}(2 + r^{2})  \nonumber  \\
 & &  -c(1 + r^{2}) - \sum_{i\not= j} 4\alpha_{i} \alpha_{j} r^{2} \}/ (1 + 
r^{2})^{2}  \ \  (1 \leq i \leq n), \label{y-i} \\
\vspace*{.1in}
\Ric(U_{j}, U_{j})& = & -\frac{f''(r)}{f(r)} + 
(2k-2)\frac{1-f'(r)^2}{f(r)^{2}} -  
\sum_{i=1}^{n}\frac{h'_{i}(r)f'(r)}{h_{i}(r)f(r)}
\nonumber  \\  & &\hspace*{2.4in} (1 \leq j \leq 2k-1). \label{u-j}
\end{eqnarray}

The Ricci tensor is positive definite since each of the diagonal terms of 
the Ricci tensor has a positive term with a factor $k$ which controls 
the negative terms~\cite{W2} and dominates the off-diagonal terms 
when $k \geq k_0(n,\alpha_i,c)$. 

\begin{rmk}\label{lie}
It was claimed in~\cite{W, W2} that
$\Ric(Y_{i}, Y_{j})=0$ if $i\neq j$. 
Actually, this is only true
when the basis $X_1,\dots ,X_n$ of the Lie algebra of $G$ satisfies
\begin{equation}\label{ijk}
\lp [X_i, X_j], [X_j, X_k] \rp = 0 \ \ \mbox{for all}\ \ i\not=k,
\end{equation}
which holds in many cases but it is unclear whether such a basis
can be always found.
Say, for the nilpotent Lie algebra of strictly upper triangular 
matrices, (\ref{ijk}) is satisfied by the so-called triangular basis
given by strictly upper triangular matrices
with $1$ at one spot and $0$'s elsewhere (see \cite{W2}).
Note that by Ado's theorem any nilpotent Lie algebra 
is isomorphic to a subalgebra of the algebra of strictly upper 
triangular matrices, but it is not always possible to choose a 
basis of the subalgebra to be a subset of the triangular basis.
Other examples of  nilpotent Lie algebras with a basis satisfying (\ref{ijk})
include all the two-step nilpotent Heisenberg algebras (see e.g.~\cite{CoG}), 
all nilpotent Lie algebra in $\dim \leq 5$~\cite{Dix}, 
and $22$ out of $24$ nilpotent Lie algebra in dimension $6$~\cite{Mor}. 
\end{rmk}

\begin{rmk}\label{ric-yi-yj}
However, all the results 
of~\cite{W} are true as stated because
$|\Ric(Y_{i}, Y_{j})|\leq  c(1+r^2)^{-1}$ which suffices
for the positivity of Ricci curvature.
The above estimate is obtained by computing 
$\Ric(Y_{i}, Y_{j})$ in the basis 
$\partial_r, U_1, \cdots U_{2k-1}, Y_1, \cdots, Y_n$.
By~\cite{W2} the terms $\lp R(Y_{i}, \partial_r )\partial_r, Y_{j}\rp$,
$\lp R(Y_{i}, U_{k})U_{k}, Y_{j}\rp$ vanish, and
\begin{eqnarray*}
\lp R(Y_{i}, Y_{k})Y_{k}, Y_{j}\rp & = &  
\frac{3}{4} h_{i}^{-1}(r) h_{k}^{-2}(r) h_{j}^{-1}(r) 
\lp [X_{i}, X_{k}], [X_{k}, X_{j}]\rp, \hspace{.1in} \mbox{if}\ i < k < j, \\
& = & \frac{1}{2} h_{i}^{-1}(r) h_{k}^{-2}(r) h_{j}^{-1}(r) \lp 
[X_{i}, X_{k}], [X_{k}, X_{j}] \rp  
\hspace{.1in} \mbox{if}\  i < j < k,\ \mbox{or}\  k < i < j.
\end{eqnarray*}
Clearly $h_{i}^{-1}(r) h_{k}^{-2}(r) h_{j}^{-1}(r)
|\lp [X_{i}, X_{k}], [X_{k}, X_{j}]\rp| \leq c (1+r^2)^{-1}$,
for some $c$ depending only on $n$ and the structure constants.
Thus, the same estimate holds for $|\Ric(Y_i,Y_j)|$.
\end{rmk}

\section{Controlling the error term}
\label{sec: error}
In this section we show that the extra terms 
$\check{\Ric}(X,Y)-{\Ric}(X,Y)$ in (\ref{ric}) 
(to which we refer to as the {\it error}) 
are controlled by some positive term from the Ricci 
curvature of the total space.

To compute the error for $\check{\Ric} (\partial_r, \partial_r)$, 
note that 
\ban
\langle T \partial_r, T \partial_r \rp & =  
&\langle T_W \partial_r, T_W \partial_r \rp 
\\
& = & \langle (\nabla_W \partial_r)^V, (\nabla_W \partial_r)^V \rp \\
& = & \langle \nabla_W \partial_r, W \rp^2 = 
\langle W, [\partial_r,W] \rp^2.
\ean
Since $[\partial_r,X_1+mV_1]$=0, it follows from (\ref{w}) that
\begin{equation} \label{bracket drw} 
[\partial_r,W]   =   -\frac{h_1 h_1' + m^2 f f'}{{h_1}^2 + m^2f^2} W.
\end{equation}
So 
\[\langle T \partial_r, T \partial_r \rp = 
\left(\frac{h_1 h_1' + m^2 f f'}
{{h_1}^2 + m^2f^2} \right)^2.\]
Also
\ban 
\langle A_{\partial_r}, A_{\partial_r} \rp & = & \langle A_{\partial_r} W, 
A_{\partial_r} W \rp = \langle (\nabla_{\partial_r} W)^H, 
(\nabla_{\partial_r} W)^H 
\rp  \\
& = & \langle \nabla_{\partial_r} W, W' \rp^2 = 
\frac 14 \langle W, [\partial_r, W'] \rp^2, 
\ean
where for the third equality we use the Koszul's formula
for the Riemannian connection 
to show that $(\nabla_{\partial_r}W)^H$ only has
$W^\prime$-component, and for the last equality we combine
the Koszul's formula with (\ref{bracket drw}). 
Now from (\ref{w'}), we compute 
\begin{equation}    \label{rw'}
[\partial_r, W']  = -\frac{{h_1} {h_1}' + m^2 f f'}{h_1^2 + m^2f^2} W' + 
\frac{f'h_1-fh_1'}{({h_1}^2 + m^2f^2)^{1/2}} \cdot 
\left(\frac{m}{{h_1}^2} X_1 + \frac{1}{f^2} V_1 \right).
\end{equation}                       
Thus $$ \langle A_{\partial_r}, A_{\partial_r} \rp  = 
\left( m \frac{f h_1' - h_1 f'}{{h_1}^2 + m^2 f^2} \right)^2. $$

%
On the other hand, since by the Koszul's formula 
$\nabla_{\partial_r}{\partial_r}=0$, we get
\ban
\langle \nabla_{\partial_r} S, \partial_r 
\rp & = & - \left( \frac{h_1 h_1' +{m}^2 f 
f'}{{h_1}^2 + {m}^2f^2} \right)' \\
& = & - 
\frac{({h_1'}^2 + h_1 h_1'' + 
{m}^2 {f'}^2 + {m}^2 ff'')
({h_1}^2 + {m}^2 f^2) -
2(h_1 {h_1}' + {m}^2 f f')^2}
{({h_1}^2 + {m}^2 f^2)^2}.
\ean
Hence, the error for $\check{\Ric} (\partial_r, \partial_r)$ is 
\ban & & 2\langle A_{\partial_r}, A_{\partial_r} \rp + 
\langle T \partial_r, T 
\partial_r \rp - \langle \nabla_{\partial_r} S, \partial_r \rp \\
& & = \frac{3{m}^2(h_1'f-f'h_1)^2}{({h_1}^2 + 
{m}^2 f^2)^2} + \frac{h_1h_1'' + {m}^2 f f''}{{h_1}^2 + 
{m}^2f^2}  \geq \frac{h_1h_1'' + {m}^2f f''}{{h_1}^2 + {m}^2f^2}  \\
& & = \frac{\frac{h_1''}{h_1} + 
{m}^2\frac{f''}{f}  (\frac{f}{h_1})^2}{1+ 
{m}^2 (\frac{f}{h_1})^2}  = 
\frac{2\alpha_1 \left[ (2\alpha_1 +1)r^2 -1 \right]-
{m}^2(\frac{3}{2}+\frac{1}{4}r^2)r^2(1+r^2)^{2^{n+1}-1}}
{\left[1+{m}^2r^2(1+r^2)^{2^{n+1}-1}\right] (1+r^2)^2}.
\ean
Here in the last equality we used some of the the following equations. 
\ba
\left(\frac{f}{h_1}\right)^2 & = & r^2(1+r^2)^{2^{n+1}-1}, \\
   \frac{h_i'}{h_i} & = & -\frac{2\alpha_i r}{1+r^2}, \ \ \ \ \ \   
\frac{h_i''}{h_i} = \frac{2\alpha_i[(2\alpha_i+1)r^2-1]}{(1+r^2)^2},  \\
   \frac{f'}{f} &  = &\frac{1+\frac{r^2}{2}}{r(1+r^2)}, \ \ \ \ \ \ 
\frac{f''}{f} = -\frac{\frac{3}{2} + \frac{1}{4}r^2}{(1+r^2)^2}.
\ea

By (\ref{p-r}), $\Ric(\partial_r, \partial_r)$ has a positive term 
$(2k-1)\frac{r^{2} + 6}{4}/(1+r^2)^2$. 
Clearly, 
\[
(2k-1)\frac{r^{2} + 6}{4} > - \frac{2\alpha_1 
\left[ (2\alpha_1 +1)r^2 -1 \right] -{m}^2(\frac{3}{2}+
\frac{1}{4}r^2)r^2(1+r^2)^{2^{n+1}-1}}{1+{m}^2r^2(1+r^2)^{2^{n+1}-1} } \]
for $k \geq \frac{2}{3}\alpha_1 + \frac{3}{2}$. 
Note that $\alpha_1 =  2^{n} -\frac{1}{4}$. 
Therefore, $\check{\Ric} (\partial_r, \partial_r)$ 
is positive for all integers $k \geq k_0 +  \frac{2^{n+1}}{3} + \frac 43$,
where $k_0$ has been chosen before Remark~\ref{lie}. 

For $\check{\Ric} (Y_i, Y_i)$, $i = 2, \cdots, n$, the error $\geq - \langle 
\nabla_{Y_i} S, Y_i \rp$. We can compute
\ban
 \langle \nabla_{Y_i} S, Y_i \rp & = & - \langle Y_i, [S, Y_i] \rp \\
 & = & -\frac{h_i'}{h_i} \cdot \frac{h_1 h_1' + 
{m}^2f f'}{h_1^2 +{m}^2 f^2} = 
-\frac{h_i'}{h_i} \cdot \frac{\frac{h_1'}{h_1}+ 
{m}^2\frac{f'}{f} (\frac{f}{h_1})^2 
}{1+ {m}^2(\frac{f}{h_1})^2}\\
& = & \frac{2\alpha_i r^2 
\left[-2\alpha_1+m^2(1+\frac{r^2}{2})(1+r^2)^{2^{n+1}-1}\right]}
{1+{m}^2r^2(1+r^2)^{2^{n+1}-1}}/ (1+r^2)^2,
\ean
which again is controlled by the positive term $(2k-1)\alpha_{i}(2 + 
r^{2})/(1+r^2)^2$ in $\Ric(Y_{i}, Y_{i})$ (see (\ref{y-i})) 
when $k \geq 1$.

For $\check{\Ric} (U_j, U_j)$, 
$j = 2, \cdots, 2k-1$, the error $\geq - \langle 
\nabla_{U_j} S, U_j \rp$. Again
\ban
\langle \nabla_{U_j} S, U_j \rp & = & -\langle U_j, [S, U_j] \rp \\
& = & -\frac{f'}{f} \cdot \frac{h_1 h_1' + {m}^2f f'}{{h_1}^2 +{m}^2 f^2}  \\
& = & - \frac{(1+\frac{r^2}{2})(-2\alpha_1+m^2
(1+\frac{r^2}{2})(1+r^2)^{2^{n+1}-1})}
{(1+{m}^2r^2(1+r^2)^{2^{n+1}-1})(1+r^2)^2}.
\ean
This is controlled by 
the positive term $(2k-2)\frac{1-{f^\prime}^2}{f^2}$ in $\Ric (U_j, U_j)$ 
(see (\ref{u-j})) when $k \geq \frac{2^{n+1}}{3} +1$. Note that 
\[
\frac{1-{f^\prime}^2}{f^2} \geq \frac{\frac 32 + r^2}{(1+r^2)^2}.
\]

For $\check{\Ric} (W',W')$, 
the error $\geq - \langle \nabla_{W'} S, W' \rp = \langle 
W', [S, W'] \rp$. 
Denote $g=\frac{h_1 h_1' + {m}^2f f'}{{h_1}^2 + {m}^2f^2}$. 
Then $S = -g\partial_r$, and $[S,W']= -g[\partial_r, W']$. 
Applying (\ref{rw'}) we have 
\ban
\langle \nabla_{W'} S, W' \rp & = 
& -g^2 +g\cdot \frac{f'h_1-fh_1'}{{h_1}^2 + {m}^2f^2} \left[ 
\frac{{m}^2f}{h_1}-\frac{h_1}{f}\right] \\
& = & -g \frac{\frac{f'}{f}+ {m}^2\frac{h_1'}{h_1} (\frac{f}{h_1})^2 
}{1+ {m}^2(\frac{f}{h_1})^2} \\
& = & 
\frac{ -\left(\frac{h_1'}{h_1} + {m}^2\frac{f'}{f} 
(\frac{f}{h_1})^2\right)  \left(\frac{f'}{f}+ {m}^2\frac{h_1'}{h_1} 
(\frac{f}{h_1})^2 \right)}{(1+ {m}^2(\frac{f}{h_1})^2)^2}. 
\ean
Now 
\ban
\lefteqn{-(1+r^2)^2\left(
\frac{h_1'}{h_1} + {m}^2\frac{f'}{f} 
\left(\frac{f}{h_1}\right)^2\right)  
\left(\frac{f'}{f}+ {m}^2\frac{h_1'}{h_1} 
\left(\frac{f}{h_1}\right)^2 \right)} \\
& & = 
-\left(-2\alpha_1+m^2(1+\frac{r^2}{2})(1+r^2)^{2^{n+1}-1} 
\right) \left(1+\frac{r^2}{2} -2\alpha_1 m^2 r^4 
(1+r^2)^{2^{n+1}-1} \right)
\\ &  & =  \left[ 2\alpha_1 m^4 r^4 
(1+\frac{r^2}{2})(1+r^2)^{2^{n+2}-2} - m^2 \left( 4\alpha_1^2 r^4 
+ (1+\frac{r^2}{2})^2 \right) (1+r^2)^{2^{n+1}-1} +  
2\alpha_1 (1+\frac{r^2}{2}) \right]
\ean
On the other hand,
\ban
\Ric (W',W') & = & \Ric \left( 
\frac{{m}fY_1-h_1U_1}{({h_1}^2 + {m}^2f^2)^{1/2}}, 
\frac{{m}fY_1-h_1U_1}{({h_1}^2 + {m}^2f^2)^{1/2}}\right) \\
& = & \frac{1}{{h_1}^2 + {m}^2f^2} \left[ {m}^2f^2 \Ric (Y_1,Y_1) + 
{h_1}^2 \Ric (U_1, U_1) 
\right] \\
& = & \frac{1}{1+ {m}^2(\frac{f}{h_1})^2} 
\left[ {m}^2\left(\frac{f}{h_1}\right)^2 \Ric (Y_1,Y_1) +  
\Ric (U_1, U_1) \right].
\ean
Now the positive term 
$\frac{1}{1+ {m}^2(\frac{f}{h_1})^2} 
(2k-1) \alpha_1 {m}^2(\frac{f}{h_1})^2 
(2 + r^{2})/(1+r^2)^2$ in 
$\Ric (W',W')$
controls $\frac{1}{(1+ {m}^2(\frac{f}{h_1})^2)^2}
\left[2\alpha_1 m^4 r^4 
(1+\frac{r^2}{2})(1+r^2)^{2^{n+2}-2} \right]/(1+r^2)^2$ when 
$k \geq 1$, and the 
positive term $\frac{1}{1+ {m}^2(\frac{f}{h_1})^2} 
(2k-2)\frac{1-{f^\prime}^2}{f^2}$ in $\Ric (W',W')$ 
controls $\frac{1}{(1+ {m}^2(\frac{f}{h_1})^2)^2} 2\alpha_1 
(1+\frac{r^2}{2})/(1+r^2)^2$ when $k \geq \frac{2^{n+1}}{3} +1$.

\section{Computing mixed terms}
\label{sec: mixed terms}
Here we show that the errors from the mixed terms either vanish,
or controlled by the positive diagonal term of the Ricci curvature of the 
total space. 

First, from (\ref{s}) one can compute that 
$[S,X]$ only has component in the $X$ 
direction for $X= \partial_r$, $Y_i$ or $U_j$. 
Also $[S,W^\prime]$ is a linear combination of $Y_1$ and $U_1$.
So the last term 
\[-\frac 12 
\left(\langle \nabla_X S, Y \rp + \langle \nabla_Y S, X \rp\right) \]
in (\ref{ric}) of the error of all mixed terms 
is zero. 

Second, note that $T_W X$ vanishes for
$X=Y_i$, $U_i$, or $W^\prime$. 
(Indeed, $T_W X=(\nabla_W X)^V$ so 
$\langle T_W X,W\rangle= \langle\nabla_W X,W\rangle=
\langle [W, X],W\rangle$, and $[W,X]=0$ since
any function of $r$ has zero derivative in the direction of $X$.)
Thus, the $T$-component of the error always vanishes
for the mixed terms, and it remains to control the
$A$-component.
 
By the remark preceding (\ref{rw'}), we have 
$A_{\partial_r} W =(\nabla_{\partial_r}W)^H= 
\langle \nabla_{\partial_r} W, W' \rp W'$.
By the Koszul's formula $\langle A_{Y_i} W, W' \rp =0$
since $[W,X]=0$ for $X=Y_i$, or $W^\prime$. 
So the error is zero and 
$\check{\Ric}( \partial_r,Y_i) = 0$ for $i=2,\cdots ,n$. 
Similarly, $\check{\Ric}( \partial_r,U_j) =0$, for $j=2,\cdots ,2k-1$. 
 
Now the error for $\check{\Ric}( Y_i,Y_j)$ with $i\ne j$, $i,j=2,\cdots,n$  
is equal to $\lp A_{Y_i}, A_{Y_j} \rp$.
By the Koszul's formula 
$\lp A_{Y_i} W,Y_k\rp =\lp W,[Y_k,Y_i]\rp/2$
and  $A_{Y_i} W$ has zero components in 
$U_j$, $\partial_r$, $W^\prime$.
Hence, $\lp A_{Y_i}, A_{Y_j} \rp$ is the sum of
the terms $\lp W,[Y_k,Y_i]\rp\lp W,[Y_k,Y_j]\rp/4$.
Since $W$ has length one, the absolute value of each of the terms 
is bounded above by the product $|[Y_i,Y_k]|\cdot |[Y_k, Y_j]|$,
which is bounded above by $c(1+r^2)^{-1}$ as in~\ref{ric-yi-yj}.
Hence, the error is controlled by the positive diagonal terms of 
the total space.
 
Similarly, by the Koszul's formula 
$\lp A_{U_i} W,U_k\rp =\lp W,[U_k,U_i]\rp/2$
and  $A_{U_i} W$ has zero components in 
$Y_j$, $\partial_r$, $W^\prime$.
Hence, $\lp A_{U_i}, A_{U_j} \rp$ is the sum of
the terms $\lp W,[U_k,U_i]\rp\lp W,[U_k,U_j]\rp/4$
which vanishes, if $i\ne j$ and $i,j=2,\cdots,2k-1$, 
because by (\ref{vbracket}) we have $[U_i, U_k]=0$. 
Thus, $\check{\Ric}(U_i,U_j)=0$ if $i\ne j$ and
$i,j=2,\cdots,2k-1$, 

By above $A_{Y_i}, A_{U_j}$
are orthogonal since they only have nonzero components
in $Y_k$'s and $U_k$'s, respectively. Thus, 
$\check{\Ric}( Y_i,U_j), \ \ i=2,\cdots,n,\ j= 2,\cdots ,2k-1$ 
is zero. 

Finally, the errors for
\ban 
& & \check{\Ric}( \partial_r,W'), \\
& & \check{\Ric}( Y_i,W'), \ \ i=2,\cdots,n,  \\
& & \check{\Ric}( U_j,W'), \ \ j=2,\cdots,2k-1
\ean
vanish since $\nabla_{W'}W$ only has a nonzero
component in $\partial_r$, 
and 
$\langle \nabla_{Y_i}W, \partial_r \rp = 
\langle  \nabla_{U_j}W, \partial_r \rp =0$.

This completes the proof that 
the manifold $(G \times {\mathbb  C}^k )/\Phi_m(\mathbb R)$ 
has positive Ricci curvature for all integers $k \geq k_0 + \frac{2^{n+1}}{3} + \frac 43$.

\begin{rmk}
Note that our $k$ depends only on $n,c, \alpha_i$, and is not optimal. 
There are other choices of $\alpha_i$ 
(or even $h_i(r)$ which work in special cases) 
that can make $k$ smaller (cf.~\cite{W2}).
\end{rmk}

\section{Metrics on vector bundles over nilmanifolds}
\label{sec: metrics on bundle}
Here is a technical version of~\ref{main thm} which we need for
the applications.

\begin{thm}\label{main thm tech}
Let $N$ be a nilmanifold and $\alpha\in H^2(N,\mathbb Z)$.
For an integer $m$, let $L_m$ be a (unique up to isomorphism) 
complex line bundle with the first Chern class $c_1(L_m)=m\alpha$.
Let $\xi_{k,m}$ be the Whitney sum of $k$ copies of $L_m$.
Then there exists $K_0=K_0(N,\alpha)$ such that for each $k\ge K_0$
the total space of $\xi_{k,m}$ admits a complete Riemannian metric 
with $\Ric>0$.
\end{thm}
\begin{proof}
By~\cite{Mal} $N$ is diffeomorphic to the product $N_c\times\mathbb R^s$ 
where $N_c$ is a compact nilmanifold. 
Let $P$ be a principal $U(1)$-bundle over 
$N_c$ with the first Chern class equal to $\alpha$.
(Recall that $c_1$ gives a bijection between the set $[X,BU(1)]$ of 
$U(1)$-bundles over a finite cell complex $X$ and the group
$H^2(X,\mathbb Z)$.)

Being a principal circle bundle over a compact nilmanifold,
$P$ is a compact nilmanifold so by~\cite{Mal} $P=G_c/\Gamma$ 
where $G_c$ is a simply-connected nilpotent Lie group, and
$\Gamma$ is a discrete subgroup of $G$. 
Let $z$ be the (central) element of $\Gamma$ corresponding to
a fiber of $P$. By~\cite[5.1.5]{CoG}, $z$ lies in the center of $G_c$.
Let $Q$ be the (necessarily central) one-parametric subgroup containing $z$.
Consider the nilpotent Lie group $G=G_c\times\mathbb R^s$
and think of $Q$, $\Gamma$ as subgroups of $G$. 
The subgroup $Q\cap\Gamma$ of $\Gamma$ is generated by $z$
(else $z$ would be a power of some element of $\Gamma$
which is then projected to a torsion element of 
the torsion free group $\pi_1(N)=\Gamma/\langle z\rangle$).
Thus, the group $\Gamma/\langle z\rangle$
is a discrete subgroup of $G/Q$ with factor space $N$.
Now $G/\Gamma$ gets the structure of a principal
$Q/\langle z\rangle\cong U(1)$-bundle over $N$
which is isomorphic to $P\times \mathbb R^s$.
(Since the bundle projections induce the same
maps of fundamental groups, the bundles are 
fiber homotopy equivalent which implies isomorphism
for principal circle bundles.)

Now let $X_1$ be a nonzero element of the Lie algebra of $G$
tangent to $Q$. Use $X_1$ to define the $\mathbb R$-action
$\Phi_m$ as in (\ref{action}). 
Thus, there exists a positive $K_0=K_0(P)$ such that for all $k\ge K_0$, 
the manifold $X=(G\times\mathbb C^k)/\Phi_m(\mathbb R)$ has a 
complete $G$-invariant metric with $\Ric>0$, so we can consider
the quotient of $X$ by the isometric action of $\Gamma\le G$.
Note that the kernel of the $G$-action is $K_m$, hence the
kernel of the $\Gamma$-action is $K_m\cap\Gamma=\langle z\rangle$.

The actions of $G$ and $\Phi_m(\mathbb R)$ 
on $G\times\mathbb C^k$ commute, so $X/\Gamma$ is
the quotient of $G/\Gamma\times\mathbb C^k$
by the drop of $\Phi_m$. 
The drop of $\Phi_m$ can be described as 
the $Q/\langle z\rangle$-action on the $G/\Gamma$-factor,
and the same diagonal $U(1)$-action on $\mathbb C^k$.

Thus, $X/\Gamma$ is the total space of a complex vector bundle over $N$.
The bundle is the Whitney sum of $k$ copies of the same line bundle
which we denote $L_m^\prime$. 
Finally, we show that the first Chern class of $L_m^\prime$ is 
$m \alpha$.
Indeed, since $L_1^\prime$ is constructed via the standard $U(1)$-action
on $\mathbb C$, we get $c_1(L_1^\prime)=\alpha$.
Look at the associated circle bundles $S(L_1^\prime)$ and $S(L_m^\prime)$.
There is a covering  $S(L_1^\prime)\to S(L_m^\prime)$ which is an $m$-fold
covering on every fiber. 
Then the classifying map for $S(L_m^\prime)$ can be obtained as
the classifying map for $S(L_1^\prime)$ followed by the selfmap 
of $BS^1$ induced by the $m$-fold selfcover of $S^1$.
By the homotopy sequence of the universal $S^1$-bundle,
this selfmap of $BS^1$ induces the multiplication by $m$
on the second homotopy group $\pi_2(BS^1)$.
By Hurewitz it induces the same map on $H_2(BS^1,\mathbb Z)$
and hence on $H^2(BS^1,\mathbb Z)$. 
Since the first Chern class of the universal bundle generates
$H^2(BS^1,\mathbb Z)$, we deduce $c_1(L_m^\prime)=m c_1(L_1^\prime)=m\alpha$.

Since the first Chern class determines (the isomorphism type of)
the complex line bundle, 
$L_m^\prime$ is isomorphic to $L_m$, 
and the proof is complete.
\end{proof}

\begin{proof}[Proof of~\ref{cor nil}]
\label{proof of cor nil}
Let $\alpha\in H^2(N,\mathbb Z)$ be a class with a nontorsion
cup-square. 
Let $L_m$ be a (unique up to isomorphism) 
complex line bundle over $N$ with $c_1(L)=m\alpha$. 
By~\ref{main thm tech} for all sufficiently large $k$
the Whitney sum $\xi_{k,m}$ of $k$ copies of
$L_m$ has a complete metric with $Ric>0$.
The total Pontrjagin class of $\xi_{k,m}$ can be computed as
\[p(\xi_{k,m})=c(\xi_{k,m}\oplus\bar{\xi}_{k,m})=c(L_m)^k c(\bar{L}_m)^k=
(1+c_1(L_m))^k(1-c_1(L_m))^k=(1-c_1(L_m)^2)^k,\]
hence the first Pontrjagin class $\xi_{k,m}$ is given by 
$p_1(\xi_{k,m})=-k c_1(L_m)^2=-k m^2 \alpha^2$.

It remains to show that the total spaces $E(\xi_{k,m})$ of
$\xi_{k,m}$ are pairwise nonhomeomorphic for different values of $|m|$. 
Being the quotient of a Lie group by a discrete subgroup,
$N$ is parallelizable (because the tangent bundle to
any Lie group has a left-invariant trivialization
which descends to quotients by discrete subgroups).
Hence, the tangent bundle $TE(\xi_{k,m})$
has first Pontrjagin class $-km^2\alpha^2$.
Assume that 
$h\co E(\xi_{k,{m}^\prime})\to E(\xi_{k,m})$ is a homeomorphism.
Since rational Pontrjagin classes are topologically invariant,
$-km^2 h^*(\alpha^2)$ is equal to $-k(m^\prime)^2\alpha^2$ up to elements
of finite order.
Thus, the projections of the classes to $H^4(N,\mathbb Z)/Tors$
are mapped to each other by the automorphism $h^*$ of
$H^4(N,\mathbb Z)/Tors$.
Since any automorphism of a finitely generated,
free abelian group preserves multiplicities of
primitive elements, $m=\pm m^\prime$, and we are done. 
\end{proof}

\begin{proof}[Proof of~\ref{prod with comp}]
Since $\dim(T)\ge 4$, we can find $\alpha\in H^2(T,\mathbb Z)$
with a nontorsion cup-square (see~\ref{sympl nil}). 
Furthermore, we can choose $\alpha$ so that $\alpha^3=0$,
say, by taking $\alpha$ to be the pullback of a class with
a nontorsion cup-square via an obvious projection of $T$ onto
the $4$-torus. 

As in proof of~\ref{cor nil}, we find a vector bundle $\xi_{k,m}$
over $T$ such that $p_1(\xi_{k,m})=-k m^2 \alpha^2$ and
$\Ric(E(\xi_{k,m}))>0$. Since $\alpha^3=0$, all the higher
Pontrjagin classes of $\xi_{k,m}$ vanish.
Since $p_1(\xi_{k,m})$ is a nontorsion class,
$\xi_{k,m}$ does not become stably trivial after passing to
a finite cover because finite covers induce injective maps on
rational cohomology. 

Now the product $C\times E(\xi_{k,m})$ has a
(product) metric with $\Ric\ge 0$ which by~\cite{BeB}
can be assumed of $\Ric>0$ after replacing $E(\xi_{k,m})$ with
$E(\xi_{k,m})\times\mathbb R^3$ (or equivalently,
after replacing $\xi_{k,m}$ by its Whitney sum with
a trivial rank $3$ bundle). 
By~\cite[1.2]{BK2} the manifold $C\times E(\xi_{k,m})$
admits no metric with $sec\ge 0$.

It remains to show that the manifolds $C\times E(\xi_{k,m})$ 
are pairwise nonhomeomorphic for different values of $|m|$.
The tangent bundle to $C\times E(\xi_{k,m})$ is isomorphic
to $TC\times TE(\xi_{k,m})$.
Let $i\ge 0$ be the largest number such that the $i$th
integral Pontrjagin class $p_i(TC)$ is nontorsion. 
Since $T$ is parallelizable, the bundle
$TE(\xi_{k,m})$ has first integral Pontrjagin class $-km^2\alpha^2$,
and zero higher Pontrjagin classes.
Up to elements of order two, 
the total Pontrjagin classes satisfy the product formula
$p(TC\times TE(\xi_{k,m}))=p(TC)\times p(TE(\xi_{k,m}))$.
In particular, $p_{i+1}(TC\times TE(\xi_{k,m}))=
p_i(TC)\times p_1(TE(\xi_{k,m}))$ is a nontorsion class
since it is a product of nontorsion classes.

Assume that 
$h\co C\times E(\xi_{k,{m}^\prime})\to C\times E(\xi_{k,m})$ 
is a homeomorphism.
Since rational Pontrjagin classes are topologically invariant, 
$h^*$ maps $p_i(TC)\times (-km^2\alpha^2)$ 
to $p_i(TC)\times (-k(m^\prime)^2\alpha^2)$ up to elements
of finite order.
Thus, the projections of the classes to $H^{4(i+1)}(C\times T,\mathbb Z)/Tors$
are mapped to each other by the automorphism $h^*$ of
$H^{4(i+1)}(N,\mathbb Z)/Tors$.
Since any automorphism of a finitely generated,
free abelian group preserves multiplicities of
primitive elements, $m=\pm m^\prime$, and we are done. 
\end{proof}

\begin{rmk}
The above proof applies verbatim
when $C$ is a complete manifold with $Ric\ge 0$
which is the total space of a vector bundle over a compact
smooth manifold. Say, we can take $C$ to be any Ricci
positively curved vector bundles of~\cite{Nash, BerB}.
\end{rmk}

\appendix
\section{Existence of line bundles over nilmanifolds}
\label{sec: append}
Let $N$ be a nilmanifold.
In this appendix we discuss when 
there is $\alpha\in H^2(N,\mathbb Z)$ with nontorsion
cup-square $\alpha^2$. As we note in~\ref{lem ina app},
this is equivalent to the existence of a complex line
bundle $L$ such that, for any $k$, the Whitney sum of 
$k$ copies of $L$ is nontrivial as a real vector bundle.
Since any nilmanifold is diffeomorphic to the product
of a compact nilmanifold and a Euclidean space, 
it suffices to discuss the case when $N$ is compact
which we assume from now on.
Clearly, one necessary condition for the existence of $\alpha$
as above is $\dim (N)\ge 4$, 
however it is not sufficient, in general, as we show in~\ref{five-nil}.

\begin{ex}\label{sympl nil}
One situation in which $\dim (N)\ge 4$ is sufficient
is when $N$ is a torus, or more generally,
a symplectic nilmanifold. 
(Indeed, since nondegeneracy is an open condition,
the symplectic form $\omega$ 
can be approximated by a nondegenerate rational form $\omega^\prime$
(where ``rational'' means that 
the cohomology class of $\omega^\prime$ lies in the
image of $H^2(N,\mathbb Q)\to H^2(N,\mathbb R)$).
Finally, we take $\alpha$ to be the cohomology class of a 
suitable integer multiple of $\omega^\prime$.)
We refer to~\cite{OT} for concrete
examples of symplectic nilmanifolds.
Of course, if $N$ carries a cohomology $2$-class with a nontorsion
cup-square, then so does any product $N\times N^\prime$. 
\end{ex}

\begin{ex}
If $N_1$, $N_2$ are arbitrary 
compact nilmanifolds of dimension $\ge 2$,
then there is $\alpha\in H^2(N_1\times N_2,\mathbb Z)$ with
nontorsion cup-square. 
To see that first note if $N$ is a compact
nilmanifold of dimension $\ge 2$, then
the second Betti number of $N$ is positive. 
(Indeed, the conclusion is obvious if $N$ is a $2$-torus, so
we can assume $\dim(N)\ge 3$.
Since $N$ can be obtained as an iterated principal circle bundle,
there is a $\pi_1$-surjective map of $N$ onto a $2$-torus.
In particular, the first Betti number of $N$ is $\ge 2$.
Since $\dim(N)\ge 3$, the second Betti number satisfies 
$b_2(N)\ge b_1(N)^2/4$~\cite{FHT} which implies $b_2(N)>0$.) 
Then given compact nilmanifolds $N_1$, $N_2$ of dimensions $\ge 2$, 
we can choose a  class $\omega_i\in H^2(N_i,\mathbb Z)$
so that there is $F_i\in H_2(N_i,\mathbb Z)$ with
$\langle \omega_i,F_i\rangle\neq 0$.
Let $p_i\co N_1\times N_2\to N_i$ be the projection.
Now one of the classes $p_1^*\omega_1$, $p_2^*\omega_2$, 
$p_1^*\omega_1+p_2^*\omega_2$
has a nontorsion cup-square because if $p_1^*\omega_1^2$, 
$p_2^*\omega_2^2$
are torsion classes, then up to elements of finite order
$(p_1^*\omega_1+p_2^*\omega_2)^2=\omega_1\times\omega_2$, and 
$\langle\omega_1\times\omega_2,F_1\times F_2\rangle=
\langle\omega_1,F_1\rangle\cdot
\langle\omega_2,F_2\rangle\neq 0.$ 
\end{ex}

\begin{coex}\label{five-nil} We now present a compact
$5$-dimensional nilmanifold $N$ such that any class in
$H^2(N,\mathbb Z)$ has a zero cup-square.
Let $x_1,x_2,x_3,x_4$ be generators of $H^1(T,\mathbb Z)$ where
$T$ is the $4$-torus. Let $p\co N\to T$ be a principal circle bundle over
$T$ with Euler class $e=x_1x_2+x_3x_4.$ Look at the Gysin
sequence of the bundle with integer coefficients:
\[
H^0(T,\mathbb Z)\stackrel{\cup e}\rightarrow
H^2(T,\mathbb Z)\stackrel{p^*}\rightarrow
H^2(N,\mathbb Z)\rightarrow 
H^1(T,\mathbb Z) \stackrel{\cup e}\rightarrow H^3(T,\mathbb Z)
\rightarrow\dots \]
By a direct computation $H^1(T,\mathbb Z) 
\stackrel{\cup e}\rightarrow H^3(T,\mathbb Z)$
is an isomorphism (this also follows from the Hard Lefschetz theorem since
$e$ is proportional to the standard K\"ahler form for $T$).
By exactness $p^*$ is onto. Now take an arbitrary
$y\in H^2(N,\mathbb Z)$ and find $x$ with $y=p^*(x)$.
Since $H^4(T,\mathbb Z)\cong\mathbb Z$, the class $x^2$ is a
rational multiple of $e^2$, hence 
$y^2=(p^*x)^2=p^*(x^2)$ is a rational multiple of
$p^*(e^2)=(p^* e)^2=0$ where the last inequality
holds because $e$ is in the kernel of $p^*$.
\end{coex}

\begin{lem}\label{lem ina app}
Given a space $X$ homotopy equivalent to a finite cell complex,
the following are equivalent:

(i) there exists $\alpha\in H^2(X,\mathbb Z)$ whose
cup-square $\alpha^2$ is nontorsion;

(ii) there exists a complex line bundle $L$ over $X$ such that,
for any $k$, the Whitney sum of $k$ copies of $L$ is
nontrivial as a real vector bundle.
\end{lem}
\begin{proof} Assume (i), and find a complex line bundle 
$L$ over $X$ with first Chern class $c_1(L)=\alpha$. 
Let $\xi_k$ be the Whitney sum of $k$ copies of $L$.
Hence $p_1(\xi_k)=-k c_1(L)^2=-k\alpha^2$ is nontorsion. 
This proves (ii).

Now assume (ii) and take $\alpha=c_1(L)$. 
Arguing by contradiction, assume $\alpha^2$ is a torsion class, 
that is $k\alpha^2=0$ for some $k>0$. 
As above we get $p_1(\xi_k)=-k\alpha^2=0$, hence all
rational Pontrjagin classes vanish.
In particular, the Pontrjagin character $ph$ of $\xi_k$ vanishes.
Recall that $ph$ defines an isomorphism of 
$\widetilde{KO}(X)\otimes\mathbb Q$ and 
$\oplus_{i>0} H^{4i}(X,\mathbb Q)$, thus $\xi_k$ has finite order
in $\widetilde{KO}(X)$. In particular, for some $m>0$,
the Whitney sum of $m$ copies of $\xi_k$ is trivial,
as a real bundle. The contradiction proves (i).
\end{proof}

\small
\bibliographystyle{amsalpha}
\bibliography{rb}

\providecommand{\bysame}{\leavevmode\hbox to3em{\hrulefill}\thinspace}
\begin{thebibliography}{FHT97}

\bibitem[And90]{And}
M.~Anderson, \emph{On the topology of complete manifolds of nonnegative {R}icci
  curvature}, Topology \textbf{29} (1990), no.~1, 41--55.

\bibitem[BB78]{BerB}
L.~B{\'e}rard-Bergery, \emph{Certains fibr\'es \`a courbure de {R}icci
  positive}, C. R. Acad. Sci. Paris S\'er. A-B \textbf{286} (1978), no.~20,
  A929--A931.

\bibitem[BB86]{BeB}
L.~B{\'e}rard-Bergery, \emph{Quelques exemples de vari{\'e}t{\'e}s
  riemanniennes compl{\`e}tes non compactes {\`a} courbure de {R}icci
  positive}, C. R. Acad. Sci. Paris S{\'e}r. I Math. \textbf{302} (1986),
  no.~4, 159--161.

\bibitem[Bes87]{Bes}
A.~L. Besse, \emph{Einstein manifolds}, Springer-Verlag, 1987.

\bibitem[BKa]{BK3}
I.~Belegradek and V.~Kapovitch, \emph{Obstructions to nonnegative curvature and
  rational homotopy theory}, 2000, submitted, available electronically at the
  xxx-achive: http://arxiv.org/abs/math.DG/0007007.

\bibitem[BKb]{BK2}
I.~Belegradek and V.~Kapovitch, \emph{Topological obstructions to nonnegative
  curvature}, 1999, to appear in Math.~Ann., available electronically at the
  xxx-achive: http://arxiv.org/abs/math.DG/0001125.

\bibitem[CG90]{CoG}
L.~J. Corwin and F.~P. Greenleaf, \emph{Representations of nilpotent {L}ie
  groups and their applications. {P}art {I}, {B}asic theory and examples},
  Cambridge University Press, 1990.

\bibitem[Dix58]{Dix}
J.~Dixmier, \emph{Sur les repr\'esentations unitaries des groupes de {L}ie
  nilpotents. {I}{I}{I}}, Canad. J. Math. \textbf{10} (1958), 321--348.

\bibitem[FHT97]{FHT}
M.~Freedman, R.~Hain, and P.~Teichner, \emph{Betti number estimates for
  nilpotent groups}, Fields Medallists' lectures, World Sci. Publishing, 1997,
  pp.~413--434.

\bibitem[Gro93]{Grl}
D.~Gromoll, \emph{Spaces of nonnegative curvature}, Differential geometry:
  Riemannian geometry, Amer. Math. Soc., 1993, pp.~337--356.

\bibitem[Mal51]{Mal}
A.~I. Mal'cev, \emph{On a class of homogeneous spaces}, Amer. Math. Soc.
  Translation \textbf{1951} (1951), no.~39, 33.

\bibitem[Mor58]{Mor}
V.~V. Morozov, \emph{Classification of nilpotent {L}ie algebras of sixth
  order}, Izv. Vys\v s. U\v cebn. Zaved. Matematika \textbf{1958} (1958), no.~4
  (5), 161--171.

\bibitem[Nas79]{Nash}
J.~C. Nash, \emph{Positive {R}icci curvature on fibre bundles}, J. Differential
  Geom. \textbf{14} (1979), no.~2, 241--254.

\bibitem[{\"O}W94]{OW}
M.~{\"O}zaydin and G.~Walschap, \emph{Vector bundles with no soul}, Proc. Amer.
  Math. Soc. \textbf{120} (1994), no.~2, 565--567.

\bibitem[TO97]{OT}
A.~Tralle and J.~Oprea, \emph{Symplectic manifolds with no {K}\"ahler
  structure}, Lecture Notes in Mathematics, 1661, Springer-Verlag, 1997.

\bibitem[Wei88]{W}
G.~Wei, \emph{Examples of complete manifolds of positive {R}icci curvature with
  nilpotent isometry groups}, Bull. Amer. Math. Soc. (N.S.) \textbf{19} (1988),
  no.~1, 311--313.

\bibitem[Wei89]{W2}
G.~Wei, \emph{Aspects of positively {R}icci curved spaces: new examples and the
  fundamental group}, Ph.D. thesis, State University of {N}ew {Y}ork at {S}tony
  {B}rook, 1989, available at http://www.math.ucsb.edu/$\sim\!\!$~wei/.

\end{thebibliography}

\

DEPARTMENT OF MATHEMATICS, 253-37, CALIFORNIA INSTITUTE OF TECHNOLOGY,
PASADENA, CA 91125, USA

{\normalsize
{\it email:} \texttt{ibeleg@its.caltech.edu}}

\

DEPARTMENT OF MATHEMATICS, UNIVERSITY OF CALIFORNIA SANTA BARBARA,
SANTA BARBARA, CA 93106, USA

{\normalsize
{\it email:} \texttt{wei@math.ucsb.edu}}

\end{document}